\documentclass[graybox]{svmult}
\usepackage{type1cm}
\usepackage{makeidx}       
\usepackage{graphicx}
\usepackage{multicol}
\usepackage[bottom]{footmisc}
\usepackage{newtxtext}       
\usepackage{newtxmath}   
\usepackage{subfigure}
\usepackage{mathtools}

\newcommand{\nono}{\nonumber}

\newcommand{\md}{\partial^\bullet_t}
\newcommand{\grad}{\nabla}
\newcommand{\Ip}{\mathbb{I}}

\newcommand{\ter}[1]{\textcolor{red}{#1}}

\newcommand{\Er}[1]{\mathcal{E}_{#1}}

\renewcommand{\nu}{\vec{n}}
\renewcommand{\langle}{(}
\renewcommand{\rangle}{)}
\newcommand{\wb}{w_b}
\newcommand{\voh}{V^h}
\newcommand{\soh}{V_0^h}
\newcommand{\tw}{w}%{\tilde{w}}
\newcommand{\vx}{\vec{x}}
\newcommand{\vxm}{|\vx_\rho|}

\newcommand{\vtau}{\vec{\tau}}
\newcommand{\vnu}{\vec{\nu}}
\newcommand{\vd}{\vec{\delta}}
\newcommand{\vX}{\vec{X}}

\newcommand{\vXmn}[1]{|\vX_\rho^{#1}|}
\newcommand{\vTau}{\vec{\mathcal{T}}}
\newcommand{\vNu}{\vec{\mathcal{N}}}
\newcommand{\ff}{F}
\newcommand{\vxi}{\vec{\xi}}

\newcommand{\dtn}{\Delta t_n}
\newcommand{\Dt}{\Delta t}

\newcommand{\drho}{\; {\rm d}\rho}

\def\epsilon{\varepsilon}

\makeindex 

\begin{document}

\title*{Finite element approximation of a system coupling curve evolution with prescribed normal contact to a fixed boundary to reaction-diffusion on the curve}
\titlerunning{Curve evolution with normal contact to a fixed boundary coupled to a PDE on the curve}
\author{Vanessa Styles \and James Van Yperen}
\institute{Vanessa Styles \at Department of Mathematics, University of Sussex, \email{v.styles@sussex.ac.uk}
\and James Van Yperen \at Department of Mathematics, University of Sussex, \email{j.vanyperen@sussex.ac.uk}}

\maketitle

\abstract{We consider a finite element approximation for a system consisting of the evolution of a curve evolving by forced curve shortening flow coupled to a reaction-diffusion equation on the evolving curve. The curve evolves inside a given domain $\Omega\subset \mathbb{R}^2$ and meets $\partial \Omega$ orthogonally. 
The scheme for the coupled system is based on the schemes derived in \cite{BDS17} and \cite{DE98}. We present numerical experiments and show the experimental order of convergence of the approximation.} 

%\abstract{Each chapter should be preceded by an abstract (no more than 200 words) that summarizes the content. The abstract will appear \textit{online} at \url{www.SpringerLink.com} and be available with unrestricted access. This allows unregistered users to read the abstract as a teaser for the complete chapter.}

\section{Introduction}
\label{sec:intro}
We consider a curve $\Gamma(t)$ evolving by forced curve shortening flow inside a given bounded domain $\Omega\subset \mathbb{R}^2$, with the forcing being a function of the solution, $w \text{ : } \Gamma(t) \rightarrow \mathbb{R}$, of a reaction-diffusion equation that holds on $\Gamma(t)$, such that 
\begin{align}
    v &= \kappa + f(w) & \text{ on } \Gamma(t), \text{  } t \in (0,T], \label{vel_eq} \\
    \md w &= w_{ss} + \kappa \, v \, w + g(v,w) & \text{ on } \Gamma(t), \text{  } t \in (0,T], \label{w_eq}
\end{align}
subject to the initial data 
$ \Gamma(0) = \Gamma_0$ and  $w(\cdot,0) = w_0 \text{ on } \Gamma_0$. 

Here $v$ and  $\kappa$ respectively denote the normal velocity and mean curvature of $\Gamma(t)$, corresponding to
the choice $\vnu$ of a unit normal, $s$ is the arclength parameter on $\Gamma(t)$ and $\md w:= w_t +  v \frac{\partial w}{\partial\nu}$ denotes the material derivative of $w$. % and $d\in \mathbb{R}_{>0}$.  
In addition we impose that the curve meets the boundary $\partial\Omega$ orthogonally. To this end we assume that $\partial \Omega$ is given by a smooth function $\ff$  such that %$\ff\in C^{2,1}(\mathbb{R}^2)$ such that
$$
\partial\Omega=\{p\in \mathbb{R}^2: F(p)=0\}~~\text{and}~~|\nabla F(p)|=1~~\forall~p\in \partial\Omega.
$$
Coupling the parametrisation of (\ref{vel_eq}), (\ref{w_eq}) that is presented in \cite{BDS17} for the setting in which $\Gamma(t)$ is a closed curve, with the formulation of (\ref{vel_eq}) presented in \cite{DE98} for the setting in which $\Gamma(t)$ meets the boundary $\partial\Omega$ orthogonally, yields the following system:
%  & (\vxm \tw)_t - (\psi \, \tw)_\rho =  \left( \frac{\tw_\rho}{\vxm} \right)_\rho + \vxm \, g(v,\tw), & \text{ for } (\rho,t) \in \Ip \times (0,T), \qquad \label{weq} \\
% & \tw_t - (\vx_t \cdot \vtau)\, \frac{1}{| \vx_{\rho} |} \,\tw_{\rho} - \frac{1}{| \vx_{\rho} |} \left( \frac{\tw_{\rho}}{| \vx_{\rho} |} \right)_\rho - \kappa \, v \, \tw = g(v,\tw), & \text{ for } (\rho,t) \in \Ip \times (0,T), \qquad \label{weq} \\
\begin{align}
    & \alpha \vx_t + (1 - \alpha) \langle \vx_t\cdot \vnu \rangle \, \vnu  = \frac{\vx_{\rho\rho}}{\vxm^2} + f(w) \, \vnu, \qquad \qquad \quad\,\, \text{ for } (\rho,t) \in \Ip \times (0,T) \label{xeq} \\
      & \tw_t - (\vx_t \cdot \vtau) \, \frac{1}{| \vx_{\rho} |} \,\tw_{\rho} - \frac{1}{| \vx_{\rho} |} \left( \frac{\tw_{\rho}}{| \vx_{\rho} |} \right)_\rho - \kappa \, v \, \tw = g(v,\tw),   \text{ for } (\rho,t) \in \Ip \times (0,T) \label{weq} \\
        &\vx(\rho,0) = \vx_0(\rho),~~~~\tw(\rho,0)= \tw_0(\rho) \qquad \qquad \qquad \qquad \, \, \, \, \, \text{ for } \rho \in \Ip \label{xic}\\
    &\tw(\rho,t) = \wb, \qquad \qquad \qquad \qquad \qquad \qquad \qquad \quad \text{ for } (\rho,t) \in \{0,1\} \times [0,T] \, \label{wbc} \\
           & F(\vx(\rho,t)) = 0, \qquad \qquad \qquad \qquad \qquad \qquad \qquad \, \, \, \text{ for } (\rho,t) \in \{0,1\} \times [0,T] \label{Feq} \\
        &\langle \vx_\rho(\rho,t)\cdot \nabla^\perp F(\vx(\rho,t) \rangle = 0, \qquad \qquad \qquad \qquad  \, \, \,\, \, \,\text{ for } (\rho,t) \in \{0,1\} \times [0,T].\label{xFp}
\end{align}
\iffalse
    \begin{eqnarray}
       & \alpha \vx_t + (1 - \alpha) \langle \vx_t\cdot \vnu \rangle \vnu  = \frac{\vx_{\rho\rho}}{\vxm^2} + f(w)\vnu,  \qquad\text{ for } (\rho,t) \in \Ip \times (0,T), &\label{xeq} \\
      % &\tw_t - (\vx_t \cdot \vtau)\, \frac{1}{| \vx_{\rho} |} \,\tw_{\rho} - \frac{d}{| \vx_{\rho} |} 
%\frac{\partial}{\partial \rho} 
%\left( \frac{\tw_{\rho}}{| \vx_{\rho} |} \right)_\rho 
- \kappa \, v \, \tw = g(v,\tw), \qquad \text{ for } (\rho,t) \in \Ip \times (0,T), &\label{weq} \\
      &   (\vxm \tw)_t - (\psi \, \tw)_\rho =  \left( \frac{\tw_\rho}{\vxm} \right)_\rho + \vxm \, g(v,\tw), \qquad \text{ for } (\rho,t) \in \Ip \times (0,T), &\label{weq} \\
        &\vx(\rho,0) = \vx_0(\rho),~~~~\tw(\rho,0)= \tw_0(\rho)  \qquad \text{ for } \rho \in \Ip,& \label{xic}\\
    &\tw(\rho,t) = \wb, \qquad\text{ for } (\rho,t) \in \{0,1\} \times [0,T], &\label{wbc} \\
           & F(\vx(\rho,t)) = 0,  \qquad\text{ for } (\rho,t) \in \{0,1\} \times [0,T], &\label{Feq} \\
        &\langle \vx_\rho(\rho,t)\cdot \nabla^\perp F(\vx(\rho,t) \rangle = 0,  \qquad\text{ for } (\rho,t) \in \{0,1\} \times [0,T].&\label{xFp}
\end{eqnarray}
\fi
Here $\alpha\in (0,1]$, $\Ip := (0,1)$, $\vx(\cdot,t):[0,1]\rightarrow \mathbb{R}^2$, $w(\rho,t) := w(\vx(\rho,t),t)$, $(\rho,t) \in [0,1] \times [0,T]$, and the unit tangent and unit normal to $\Gamma(t)$ are respectively given by 
%Before deriving a weak formulation of the coupled system we first note that in the parametric setting we are working in we have 
$%\begin{equation*}  \label{tau}
%\displaystyle 
\vtau = \vx_s = \frac{\vx_\rho}{|\vx_\rho|}$ and $\vnu = \vtau^\perp 
$ %\end{equation*}
%are a  unit tangent and unit normal to $\Gamma(t)$ respectively, 
where $(\cdot)^\perp$ denotes counter-clockwise rotation by $\frac{\pi}{2}$.

%\end{subequations}

The formulation of curve shortening flow in the form of \eqref{xeq} for a closed curve in $\mathbb{R}^2$ was presented and analysed in \cite{EF16}, where the DeTurck trick is used in coupling the motion of the curve to the harmonic map heat flow, with the parameter $\alpha\in (0,1]$ being such that $1/\alpha$ corresponds to the diffusion coefficient in the harmonic map heat flow. Setting $\alpha\in (0,1]$ introduces a tangential part in the velocity which, at the numerical level, gives rise to a good distribution of the mesh points along the curve. Setting $\alpha =1$ one recovers the formulation introduced and analysed in \cite{DD95}, while formally setting $\alpha = 0$ yields the approach introduced in \cite{BGN11}. In \cite{DES01} the authors derive finite element approximations of a simplified version of the  parametric coupled system \eqref{xeq}-\eqref{xFp}, and two related models. In particular,  the evolution law for the parametric system derived in \cite{DES01}, can be obtained from \eqref{xeq} by setting  $\alpha=1$, 
%the boundary $\partial \Omega$  can be expressed by setting 
$\ff(x,y)=|x|-1$, and considering a slightly different formulation of the reaction-diffusion equation \eqref{weq}. In \cite{BDS17}  the authors prove optimal error bounds for a fully discrete finite element approximation of the coupled system \eqref{xeq}-\eqref{xic} for the case where  $\Gamma(t)$ is a closed curve in $\mathbb{R}^2$. While in \cite{PS15} optimal error bounds are presented for a semi-discrete finite element approximation of an alternative formulation, which is introduced and analysed in \cite{D94}, of the coupled system \eqref{xeq}-\eqref{xic}, for the case where  $\Gamma(t)$ is a closed curve in $\mathbb{R}^2$. 
Setting $\alpha=1$ and $f(w)=0$ in \eqref{xeq} and coupling the resulting equation to \eqref{Feq}, \eqref{xFp} gives rise to the model presented and analysed in \cite{DE98}, in which optimal order error bounds for a semi-discrete finite element approximation of curve shortening flow with a prescribed normal  contact to a fixed boundary are presented. In \cite{BGN07b} the authors propose parametric finite element approximations of combined second and fourth order geometric evolution equations for curves that are connected via triple or quadruple junctions or that intersect external boundaries. %The closed curve formulation of \eqref{xeq}-\eqref{xic} that was studied in \cite{BDS17} and the formulation that we consider here in which the curve meets a fixed boundary orthogonally, can both be used to model diffusion induced grain boundary motion, \cite{H88}. This phenomenon can be observed if a polycrystalline film of metal is placed in a vapour containing another metal: atoms from the vapour diffuse into the film along the grain boundaries that separate the crystals in film, this gives rise to variations of elastic energy in the film that cause the grain boundaries to move. 
\section{Weak formulation and finite element approximation}
\label{sec:wffem}
For a weak formulation of \eqref{xeq} we multiply it by $\vxm^2\vxi$, where $\vxi \in [H^1(\Ip)]^2$ is a test function, integrate in space, use integration by parts and \eqref{xFp} to obtain 
\begin{align}
    \label{acsfwf}
    & \left( \vxm^2 \left[ \alpha \, \vx_t + (1 - \alpha) \langle \vx_t\cdot\vnu \rangle \, \vnu \, \right], \vxi \right) + \left( \vx_\rho, \vxi_\rho \right) \nono \\
    & \quad = \left[ \langle \vx_\rho\cdot \grad F(\vx) \rangle \langle \vxi\cdot \grad F(\vx) \rangle \right]_0^1 + \left( \vxm^2 f(\tw) \, \vnu, \vxi \right) \qquad \forall \, \vxi \in [H^1(\Ip)]^2,
\end{align}
where $(\cdot,\cdot)$ denotes the standard $L^2(\Ip)$ inner product. For a weak formulation of \eqref{weq}  we multiply it by $\vxm \, \eta$, where $\eta \in H_0^1(\Ip)$ is a time-independent test function, integrate in space, use integration by parts and note that $\vtau_{\rho} = \kappa \,\vnu \,|\vx_\rho|$ to obtain 
\begin{align}
    \label{awwf}
    \frac{d}{dt} \left( \vxm \, w, \eta \right) + \left( \psi \, w, \eta_\rho \right) +  \left( \frac{w_\rho}{\vxm}, \eta_\rho \right) = \left( \vxm \, g(v,w), \eta \right) \quad \forall \, \eta \in H_0^1(\Ip). 
\end{align}
Here $\psi$ is the tangential velocity of $\Gamma(t)$, such that the normal and tangential velocities of $\Gamma(t)$ are given by
$%\begin{equation*} \label{vpsi}
v= \vx_t \cdot \vnu$ and $\psi= \vx_t \cdot \vtau$. 
%\end{equation*}
We now introduce a finite element approximation of \eqref{acsfwf}, \eqref{awwf}. We first let $0 = t_0 < t_1 < \cdots < t_{N-1} < t_N = T$ be a partition of $[0,T]$ with $\Delta t_n := t_n - t_{n-1}$. Next we partition the interval $\Ip$ such that $\Ip = \cup_{j=1}^J \overline{\sigma_j}$, where $\sigma_j = (\rho_{j-1},\rho_j)$, with $h_j = \rho_j - \rho_{j-1}$. 
%We define our finite element spaces as
We set 
\begin{align*}
    %\label{Vh}
    \voh &:= \{ \chi \in C(\Ip) \text{ : } \chi_{|_{\sigma_j}} \text{ is affine, } j = 1,\dots,J \} \subset H^1(\Ip) \\
    \soh &:= \{ \chi \in V^h \text{ : } \chi(\rho_j) = 0, \text{ for } j \in \{0,J\} \} %\label{Sh}
\end{align*}
and denote the standard Lagrange interpolation operator by $I^h \text{ : } C(\Ip) \rightarrow \voh$, where $(I^h \eta)(\rho_j) = \eta(\rho_j)$, for $j=1,\dots,J$. We define the discrete inner product $(\eta_1, \eta_2)^h$ by
\begin{align*}
    %\label{innerh}
    \left( \eta_1, \eta_2 \right)^h &:= \sum_{j=1}^J \int_{\sigma_j} I^h_j(\eta_1\,\eta_2) \, \drho,
\end{align*} 
where 
%by definition $\eta_i$ are at least piecewise continuous functions of the defined partition of $\Ip$ and 
$I^h_j = I^h_{|_{\sigma_j}}$ is the local interpolation operator. Our finite element approximation of \eqref{acsfwf}, \eqref{awwf} then takes the form: 

Given $(\vX^{n-1},W^{n-1}-\wb) \in [\voh]^2\times \soh$, find $(\vX^n,W^{n}-\wb)  \in [\voh]^2 \times \soh$ such that for all 
$(\vxi^h,\eta^h) \in [\voh]^2\times  \soh$ we have
\begin{eqnarray}
    \SwapAboveDisplaySkip
     \left( \vXmn{n-1}^2 \left[ \alpha D_t \vX^n + (1 - \alpha) \langle D_t \vX^n\cdot\vNu^{n-1} \rangle \, \vNu^{n-1} \right], \vxi^h \right)^h + \left( \vX^n_\rho, \vxi^h_\rho \right)\hspace{10mm}\nono \\
    = \left[ \langle \vX^{n}_\rho\cdot\nabla F(\vX^{n})\rangle \langle \vxi^h\cdot \nabla F(\vX^{n}) \rangle \right]^1_0%\nono \\
    %& \qquad \qquad 
    + \left(\vXmn{n-1}^2 f(W^{n-1}) \, \vNu^{n-1}, \vxi^h \right)^h    \label{xfea}\\
     D_t \left[ \left( \vXmn{n} \, W^n, \eta^h \right)^h \right] +  \left( \frac{W^n_\rho}{\vXmn{n}}, \eta^h_\rho \right) + \left( \Psi^{n} \, W^{n}, \eta^h_\rho \right)^h \hspace{28mm} \nono \\
 %\qquad \qquad \qquad 
   = \left( \vXmn{n} \, g(V^n,W^{n-1}), \eta^h \right)^h~~~~~~%~ \forall \, \eta^h \in \soh, 
     \label{wfea}
\end{eqnarray}
with the additional boundary constraint
\begin{align}
    \label{dis_bc}
    F(\vX^n_{0}) =  F(\vX^n_{J}) = 0.% \qquad \text{ for } \rho \in \{0,1\}.
\end{align}
Here and in what follows we set  $D_t(a^n):=(a^{n}-a^{n-1})/\Delta t_n$ and on $\sigma_j$, $j=1,\cdots,J$ we set $\vTau^n = \frac{\vX^n_\rho}{\vXmn{n}}$ , $\vNu^n = (\vTau^n)^\perp$, $\Psi^n = D_t \vX^n \cdot \vTau^n$ and $V^n = D_t \vX^n \cdot \vNu^n$.
%Our finite element approximation of \eqref{awwf} takes the following form. Given $\vX^{n-1}, \vX^n \in [\voh]^2$ and $W^{n-1} \in \voh$, find $W^n \in \voh$ such that for $\eta^h \in \soh$ we have
%\begin{align}
%    \label{wfea}
%    & D_t \left[ \left( \vXmn{n} W^n, \eta^h \right)^h \right] + d \left( \frac{W^n_\rho}{\vXmn{n}}, \eta^h_\rho \right) + \left( \Psi^{n} W^{n}, \eta^h_\rho \right)^h \nono \\
%    & \qquad \qquad \qquad = \left( \vXmn{n} g(V^n,W^{n-1}), \eta^h \right)^h \qquad \forall \, \eta^h \in \soh,
%\end{align}
%whereby we set $W^n = g_d$ for $\rho \in \{0,1\}$.
\begin{remark}
In \cite{DE98}, rather than use the nonlinear scheme presented above to approximate (\ref{xeq}), \eqref{Feq}, \eqref{xFp}, the authors present a linear scheme in which  \eqref{Feq} is not necessarily satisfied but is instead approximated through the relation $0 = \frac{d}{dt} F(\vx) = \vx_t\cdot \grad F(\vx)$. 
\label{rem}
\end{remark}

\section{Numerical results}
\label{sec:numerics}
\subsection{Solution of the discrete system (\ref{xfea}), (\ref{dis_bc})}
We solve the resulting system of nonlinear algebraic equations arising 
at each time level from the approximation (\ref{xfea}), (\ref{dis_bc}), with  $\vxi^h = \chi_j$, $ j = 1,\dots,J-1$, $\vxi^h = \nabla^\perp F(\vX^{n}) \chi_0$ and $\vxi^h = \nabla^\perp F(\vX^{n}) \chi_J$, using the following Newton scheme, where for ease of presentation we set $\alpha =1$ and $f(w) = 0$: 

Given $\vX^{n,i-1}$, with $\vX^{n,0} = \vX^{n-1}$, we set $\vX^{n,i} := \vX^{n,i-1} + \vd^i$ such that for $ j = 1,\dots,J-1$, $\vd^i$ solves 
\begin{subequations}
\begin{align}
    \SwapAboveDisplaySkip
    & \frac{1}{2} \left( q_j^{n-1} + q_{j-1}^{n-1} \right) \frac{\vd_j^i}{\dtn} - \left( \vd_{j-1}^i - 2 \vd_j^i + \vd_{j+1}^i \right) 
    %\qquad \quad \text{ for } j = 1,\dots,J-1 
    \nono \\
    & \quad = - \frac{1}{2} \left( q_j^{n-1} + q_{j-1}^{n-1} \right) D_t(\vX^{n,i-1}_j)  + \left(\vX^{n,i-1}_{j-1} - 2 \vX^{n,i-1}_j + \vX^{n,i-1}_{j+1} \right), \label{ns} \\
%    &\frac{1}{2} \, {Q^{n-1}_j} \left[ \left\langle \frac{\vd_j^i}{\dtn}, \nabla^\perp F(\vX^{n,i-1}_j) \right\rangle + \left\langle \frac{\vX^{n,i-1}_j - \vX^{n-1}_j}{\dtn}, D \nabla^\perp F(\vX^{n,i-1}_j) \, \vd_j^i \right\rangle \right] \nono \\
%    & \quad - \left\langle (\vd_j^i)', \nabla^\perp F(\vX^{n,i-1}_j) \right\rangle - \left\langle (\vX^{n,i-1}_j)', D \nabla^\perp F(\vX^{n,i-1}_j) \, \vd_j^i \right\rangle \quad \text{ for } j = 0,J \nono \\
%    & \quad = -\frac{1}{2} \, Q^{n-1}_j \left\langle \frac{\vX^{n,i-1}_j - \vX^{n-1}_j}{\dtn}, \nabla^\perp F(\vX^{n,i-1}_j) \right\rangle + \left\langle (\vX^{n,i-1}_j)', \nabla^\perp F(\vX^{n,i-1}_j) \right\rangle \label{bc1} \\
   &\frac{1}{2} \, {q^{n-1}_0} \left[ \left\langle \frac{\vd_0^i}{\dtn}\cdot\nabla^\perp F(\vX^{n,i-1}_0) \right\rangle + \left\langle D_t(\vX^{n,i-1}_0)\cdot D_\perp^2 F(\vX^{n,i-1}_0) \, \vd_0^i \right\rangle \right] \label{bc1} \\
    & \quad - \left\langle (\vd_1^i-\vd_0^i)\cdot \nabla^\perp F(\vX^{n,i-1}_0) \right\rangle - \left\langle (\vX^{n,i-1}_1-\vX^{n,i-1}_0) \cdot D_\perp^2 F(\vX^{n,i-1}_0) \, \vd_0^i \right\rangle \nono \\
    & \quad = -\frac{1}{2} \, q^{n-1}_0 \left\langle D_t(\vX^{n,i-1}_0) \cdot \nabla^\perp F(\vX^{n,i-1}_0) \right\rangle + \left\langle (\vX^{n,i-1}_1-\vX^{n,i-1}_0)\cdot \nabla^\perp F(\vX^{n,i-1}_0) \right\rangle, \nono \\
       &\frac{1}{2} \, {q^{n-1}_{J-1}} \left[ \left\langle \frac{\vd_J^i}{\dtn}\cdot\nabla^\perp F(\vX^{n,i-1}_J) \right\rangle + \left\langle   D_t(\vX^{n,i-1}_J) \cdot D_\perp^2 F(\vX^{n,i-1}_J)  \, \vd_J^i \right\rangle \right] \label{bc1a} \\
    & \quad - \left\langle (\vd_{J-1}^i-\vd_J^i)\cdot \nabla^\perp F(\vX^{n,i-1}_J) \right\rangle - \left\langle (\vX^{n,i-1}_{J-1}-\vX^{n,i-1}_J) \cdot D_\perp^2F(\vX^{n,i-1}_J) \, \vd_J^i \right\rangle \nono \\
    & \quad = -\frac{1}{2} \, q^{n-1}_{J-1} \left\langle D_t(\vX^{n,i-1}_J) \cdot \nabla^\perp F(\vX^{n,i-1}_J) \right\rangle + \left\langle (\vX^{n,i-1}_{J-1}-\vX^{n,i-1}_J)\cdot \nabla^\perp F(\vX^{n,i-1}_J) \right\rangle, \nono \\
    &\left\langle \nabla F(\vX^{n,i-1}_0) \cdot \vd^i_0\right\rangle = -F(\vX^{n,i-1}_0) ~~\text{and}~~ \left\langle \nabla F(\vX^{n,i-1}_J)\cdot \vd^i_J \right\rangle = -F(\vX^{n,i-1}_J),
    %\qquad \qquad \qquad \qquad \qquad \text{ for } j = 0,J. 
    \label{bc2} 
\end{align}
\end{subequations}
where $q_j^{n-1} = |\vX^{n-1}_{j+1} - \vX^{n-1}_j|^2$, $ D^2_\perp = \left(
\begin{matrix} 
-\partial^2_{xy} & -\partial^2_{yy} \\
\partial^2_{xx} & \partial^2_{xy}
\end{matrix}
\right) $, and in an abuse of notation we have redefined $D_t$ from the previous section such that $D_t(\vX^{n,i-1}_j):=(\vX^{n,i-1}_j - \vX^{n-1}_j)/\dtn$. We adopt the stopping criteria $\max\limits_{j=0,J} |F(\vX_j^{n,i})| \leq \tau$ for some predetermined tolerance, $\tau$.

\subsection{Experimental order of convergence of (\ref{xfea}), (\ref{dis_bc})}

We investigate the experimental order of convergence of (\ref{xfea}), (\ref{dis_bc}) by monitoring the following errors:
\begin{align*}
     \Er{1} := \sup_{n = 0,\dots,N} \|I^h(\vx_\rho^n) - \vX_\rho^n\|_{[L^2(\Ip)]^2}^2, \quad \Er{2} := \sum_{n=1}^N \dtn \|D_t (I^h(\vx^n) - \vX^n)\|_{[L^2(\Ip)]^2}^2. %, \\
  %  & \Er{3} := \sup_{n = 0,\dots,N} \sup_{j=0,J} |F(\vX^n_j)|,
\end{align*}
In addition we show how the choice of $\alpha$ affects the size of the errors. In all examples we use a uniform mesh size $h J = 1$ and a uniform time step size $\Dt = h^2$. \\

\begin{figure}[ht!]
\centering
\subfigure{\includegraphics[width = 0.48\textwidth]{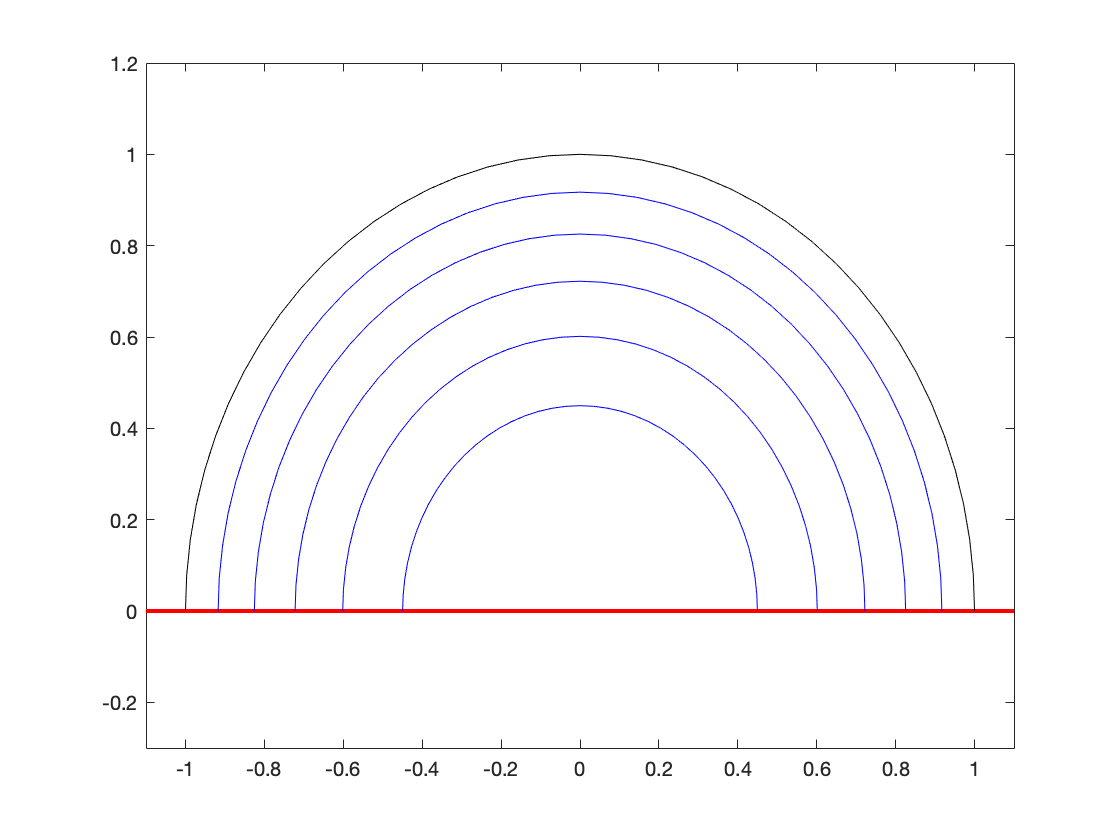}} \hspace{2mm}
\subfigure{\includegraphics[width = 0.48\textwidth]{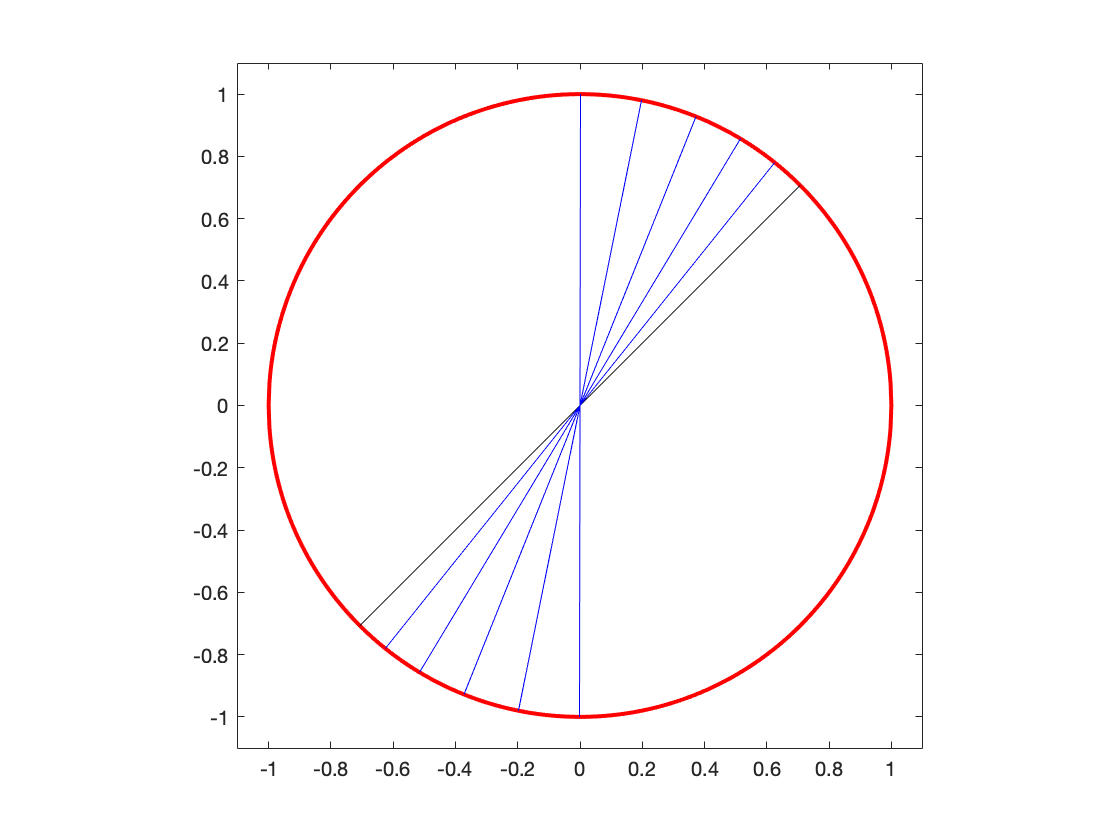}}  
    \caption{$\vec X^n$ at $t^n=0,0.08,0.16,0.24,0.32,0.4$ for the shrinking semi-circle (left), and $t^n=0,0.1,0.2,0.3,0.4,0.5$ for the rotating diameter (right).}
\label{f:ex1}
\end{figure}

%and quantifying them using the estimated order of convergence (eoc)
%\begin{align*}
%    EOC_{i,j} := \frac{\ln(\Er{i,j+1}) - \ln(\Er{i,j})}{\ln(h_{j+1}) - \ln(h_j)}, \qquad EOC_{3,j} := \frac{\ln(\Er{3,j+1}) - \ln(\Er{3,j})}{\ln(\Dt_{j+1}) - \ln(\Dt_j)}
%\end{align*}
%where $i$ corresponds to the error measurement, $j$ corresponds to the relative mesh size $h_j$ and $\Er{i,j}$ corresponds to the error of the error measurement $i$ at the $j$ level. 
{\bf Example 1: } In the first example we set $T=0.4$ and $\Omega := \mathbb{R} \times \mathbb{R}_{>0}$, such that $\partial\Omega$ is given by $F(x,y) = y$. Taking $\Gamma(0)$ to be a semi circle with radius 1, the explicit solution is given by 
\begin{align*}
    \vx(\rho,t) = \sqrt{1 - 2t} \, (\cos(\pi \rho), \sin(\pi \rho))^T.
\end{align*}
%We consider a uniform mesh size $h J = 1$ and take a uniform time step size $\Dt=h^2$, with $T = 0.4$. 
In the left-hand plot in Figure \ref{f:ex1} we display:  $\vX^0$ in black, $\vec X^n$ at $t^n = 0.08 k$, $k=1,\dots,5$, in blue, and $\partial \Omega$ in red, while in Table \ref{tab:1} we display the values of $\mathcal{E}_i$, 
$i=1,2$, for $\alpha = 1$ (left) and $\alpha = 0.5$ (right). For both values of $\alpha$ we see eocs close to four, however we note that the errors for $\alpha = 0.5$ are significantly smaller than those for $\alpha = 1$. 
\begin{table}[ht!]
\caption{Errors and eocs for the shrinking semi circle with $\alpha = 1$ (left) and $\alpha = 0.5$ (right).}%for \eqref{ns}-\eqref{bc2}, 
\label{tab:1}
    \begin{tabular}{p{0.5cm}p{0.7cm}p{1.2cm}p{0.7cm}p{1.2cm}p{0.6cm}}
    \hline\noalign{\smallskip}
    $J$ & $N$ &  $\Er{1}\times e^{-3}$ & $eoc_1$ & $\Er{2}\times e^{-4}$ & $eoc_2$  \\
    \noalign{\smallskip}\svhline\noalign{\smallskip}
    10 & 40 &  4.672 & - & 20.16 &  -  \\
    20 & 160  & 0.3997 & \ter{3.55} & 1.859 & \ter{3.44}  \\
    40 & 640  & 0.02726 & \ter{3.87} & 0.1298 & \ter{3.84} \\
    80 & 2560  & 0.001742 & \ter{3.97} & 0.008347 & \ter{3.96}  \\
    \noalign{\smallskip}\svhline\noalign{\smallskip}
    %160 & 10240  & 0.0001095 & \ter{3.99} & 0.0005254 & \ter{3.99} \\
    %\hline 
    \end{tabular}
    ~~~~~
    \begin{tabular}{p{0.5cm}p{0.7cm}p{1.3cm}p{0.7cm}p{1.2cm}p{0.6cm}}
    \hline\noalign{\smallskip}
    $J$ & $N$ &  $\Er{1}\times e^{-3}$ & $eoc_1$ & $\Er{2}\times e^{-4}$ & $eoc_2$  \\ 
    \noalign{\smallskip}\svhline\noalign{\smallskip}
    10 & 40 &  1.589 &   - & 8.884 & - \\
    20 & 160  &0.1389 & \ter{3.52} & 0.8302 & \ter{3.42}  \\
    40 & 640  &  0.009514 & \ter{3.87} & 0.05798 & \ter{3.84} \\
    80 & 2560  &  0.0006087 & \ter{3.97} & 0.003729 & \ter{3.96} \\
    \noalign{\smallskip}\svhline\noalign{\smallskip}
    %160 & 10240  & 0.00003827 & \ter{3.99} & 0.0002348 & \ter{3.99} \\
    %\hline 
    \end{tabular}
\end{table}

{\bf Example 2: } In the second example we set $T=0.5$ and $\Omega$ to be the unit disc with centre $(0,0)$, such that $\partial\Omega$ is given by $F(x,y) = \frac{1}{2}\left(x^2 + y^2 - 1\right)$. In contrast to the previous example this example has been constructed so that $| \grad F(p) | = 1$ is only satisfied on $\partial\Omega$. 
By setting $f(\rho,t) = \frac{4(\rho - \frac{1}{2})}{(1-2t)^2 + 1}$ the explicit solution is given by 
\begin{align*}
    \vx(\rho,t) = \frac{2(\rho - \frac{1}{2})}{\sqrt{(1-2t)^2 + 1}}\left(1-2t,1\right)^T,
\end{align*}
such that $\Gamma(t)$ is a rotating straight line that spans the diameter of $\Omega$. %, see the right-hand plot in Figure \ref{f:ex1}. 
%We use a uniform mesh size $h J = 1$ and a uniform time step size $\Dt = 0.1 h$, with $T = 0.5$. 
In the right-hand plot of Figure \ref{f:ex1} we display: $\vX^0$ in black, $\vec X^n$ at $t^n = 0.1k$, $k=1,\dots,5$, in blue, and $\partial \Omega$ in red, while Table \ref{tab:2} displays the errors $\mathcal{E}_i$, 
$i=1,2$, for $\alpha = 1$ (left) and $\alpha = 0.5$ (right). As in Example 1, both values of $\alpha$ exhibit eocs close to four, with the errors obtained using $\alpha=0.5$ being smaller than those obtain using $\alpha=1$. However the difference in the errors for the two values of $\alpha$ in this example is much smaller than the difference in the errors for the two values of $\alpha$ in Example 1, we believe that this is due to the fact that in this example $\vx$ is a linear function. 

To demonstrate Remark \ref{rem} in Section \ref{sec:wffem}, we include Table \ref{tab:ex2de} in which we display errors obtained using the scheme in \cite{DE98}. In particular we display $\Er{i}$, $i=1,2,3$, with 
\begin{align*}
    \Er{3} := \sup_{n=0,\dots,N} \sup_{j=0,\cdots,J} |F(\vX^n_j)|.
\end{align*}
%We don't demonstrate $\Er{3}$ for the Newton scheme because the solver measures $\Er{3}$ and stops at machine tolerance after a few iterations. 
%We did not include such a table in Example 1,  since in that example $|\grad F| \equiv 1$.%
%as the scheme in \cite{DE98} obtained machine tolerance on the boundary since in this example $|\grad F| \equiv 1$. % enables the boundary conditions \eqref{bc1}-\eqref{bc2} to decouple and solve explicitly for $y=0$. 

%To this end in Table \ref{tab:ex2de} we display the errors $\Er{i}$, $i=1,2,3$, for the scheme proposed in \cite{DE98}. 
Comparing the errors in Tables \ref{tab:2} and \ref{tab:ex2de} we see that the magnitude of the errors for the Newton scheme, \eqref{ns}-\eqref{bc2}, are significantly smaller than the errors for the linear scheme in \cite{DE98}. 
%We speculate that this is due to \eqref{dis_bc} being approximated by \eqref{bc2}, rather than by linearising, and not down to the fact that we are using a Newton solver since the linear system gains very similar errors for example one. 
%We consider an example whereby $|\grad F(p)| = 1$ only for $p \in \partial \Omega$. Taking $f(\rho,t) = \frac{4(\rho - \frac{1}{2})}{(1-2t)^2 + 1}$ we see that the explicit solution is given by
%\begin{align*}
%    \vx(\rho,t) = \frac{2(\rho - \frac{1}{2})}{\sqrt{(1-2t)^2 + 1}}\left(1-2t,1\right)^T.
%\end{align*}
\begin{table}[t!]
\caption{Errors and eocs for the rotating diameter with $\alpha = 1$ (left) and $\alpha = 0.5$ (right).}%for \eqref{ns}-\eqref{bc2}, }
\label{tab:2}
    \begin{tabular}{p{0.5cm}p{0.7cm}p{1.3cm}p{0.7cm}p{1.3cm}p{0.6cm}}
    \hline\noalign{\smallskip}
    $J$ & $N$ &  $\Er{1}\times e^{-4}$ & $eoc_1$ & $\Er{2}\times e^{-5}$ & $eoc_2$ \\
    \noalign{\smallskip}\svhline\noalign{\smallskip} 
    10 & 50 & 1.440 & -&  3.040 & -  \\
    20 & 200  & 0.09198 & \ter{3.97} & 0.1925 & \ter{3.98} \\
    40 & 800 & 0.005780 & \ter{3.99} & 0.01207 & \ter{4.00}  \\
    80 & 3200 & 0.0003617 & \ter{4.00} & 0.0007552 & \ter{4.00} \\
    \noalign{\smallskip}\svhline\noalign{\smallskip}
    %160 & 800 & 0.005790 & \ter{2.00} & 0.001208 & \ter{2.00} \\
    %\hline 
    \end{tabular}
    ~~~
    \begin{tabular}{p{0.5cm}p{0.7cm}p{1.3cm}p{0.7cm}p{1.3cm}p{0.6cm}}
    \hline\noalign{\smallskip}
    $J$ & $N$ & $\Er{1}\times e^{-4}$ & $eoc_1$ & $\Er{2}\times e^{-5}$ & $eoc_2$ \\
    \noalign{\smallskip}\svhline\noalign{\smallskip} 
    10 & 50  & 1.181 & - & 2.710 & - \\
    20 & 200 & 0.07459 & \ter{3.98} & 0.1716 & \ter{3.98} \\
    40 & 800 &  0.004674 & \ter{4.00} & 0.01076 & \ter{4.00} \\
    80 & 3200 & 0.0002923 & \ter{4.00} & 0.0006727 & \ter{4.00} \\
    \noalign{\smallskip}\svhline\noalign{\smallskip} 
    %160 & 800 & 0.004678 & \ter{2.00} & 0.01076 & \ter{2.00} \\
    %\hline 
    \end{tabular}
\end{table}
\begin{table}[ht!]
\caption{Errors and eocs for the rotating diameter for scheme proposed in \cite{DE98}.}
\label{tab:ex2de}
    \begin{tabular}{p{0.5cm}p{0.7cm}p{1.1cm}p{0.7cm}p{1.1cm}p{0.7cm}p{1.0cm}p{0.6cm}}
    \hline\noalign{\smallskip}
    $J$ & $N$ & $\Er{1} \times e^{-4}$ & $eoc_1$ & $\Er{2} \times e^{-5}$ & $eoc_2$ & $\Er{3} \times e^{-3}$ & $eoc_3$ \\
    \noalign{\smallskip}\svhline\noalign{\smallskip} 
    10 & 50 & 43.83 &  & 76.20 &  & 5.771 &  \\
    20 & 200 & 3.175 & \ter{3.79} & 5.442 & \ter{3.81} & 1.563 & \ter{0.94} \\
    40 & 800 & 0.2076 & \ter{3.93} & 0.3542 & \ter{3.94} & 0.3989 & \ter{0.99} \\
    80 & 3200 & 0.01317 & \ter{3.98} & 0.02243 & \ter{3.98} & 0.1003 & \ter{1.00} \\
    \noalign{\smallskip}\svhline\noalign{\smallskip} 
    %160 & 800 & 0.02104 & \ter{1.98} & 0.03577 & \ter{1.98} & 0.3989 & \ter{0.99} \\
    %\hline 
    \end{tabular}
\end{table}

\subsection{Experimental order of convergence of the coupled scheme (\ref{xfea})-(\ref{dis_bc})}

\begin{figure}[ht!]
\centering
\subfigure{\includegraphics[width = 0.48\textwidth]{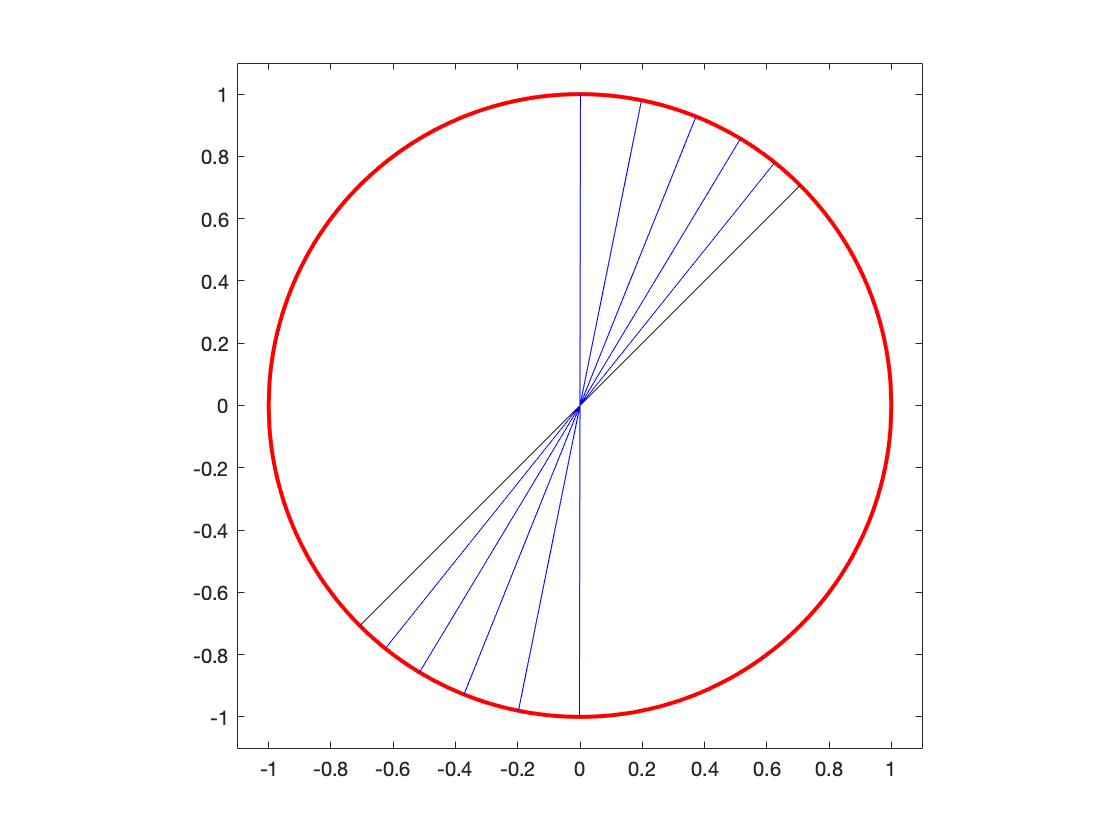}} \hspace{2mm}
\subfigure{\includegraphics[width = 0.48\textwidth]{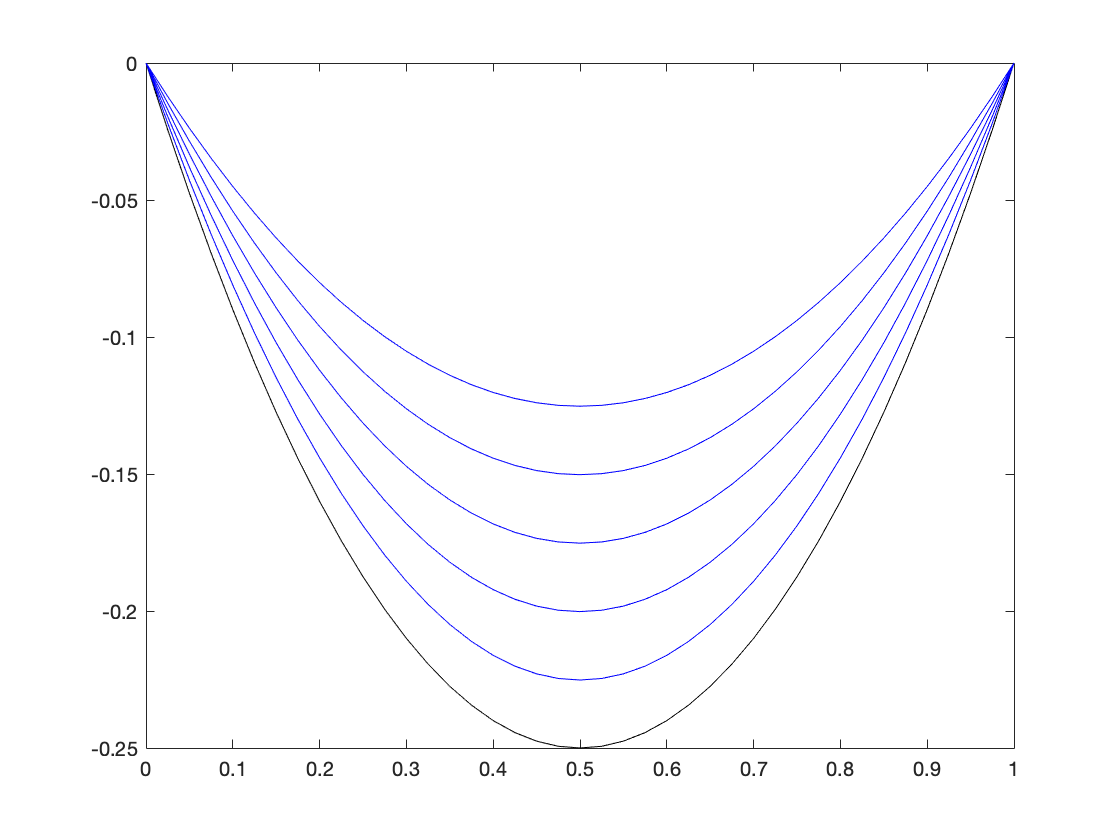}}  
\caption{$\vec X^n$ at $t^n=0,0.1,0.2,0.3,0.4,0.5$ for the rotating diameter (left) and $W^n$ at $t^n=0,0.1,0.2,0.3,0.4,0.5$ for the shrinking parabolic (right).}
\label{f:ex2}
\end{figure}

We conclude our numerical results by investigating the experimental order of convergence of the coupled scheme \eqref{xfea}-\eqref{dis_bc}. In addition to monitoring the errors $\Er{i}$, $i=1,2$, we also monitor
\begin{align*}
    \Er{4} := \sup_{n=0,\dots,N} \|I^h(w^n) - W^n\|_{L^2(\Ip)}^2, \quad \Er{5} := \sum_{n=1}^n \dtn \|I^h(w^n_\rho) - W^n_\rho\|_{L^2(\Ip)}^2.
\end{align*}
We adopt the same set-up as in Example 2, with $T=0.5$ and $\Omega$ being the unit disc with centre $(0,0)$, such that $\partial\Omega$ is given by $F(x,y) = \frac{1}{2}\left(x^2 + y^2 - 1\right)$. Setting $f(\rho,t) = \frac{4 \left( \rho^2 - \frac{w(\rho,t)}{1 - t} - \frac{1}{2} \right)}{(1-2t)^2 + 1} $ and $g = \frac{t-1}{2} - \frac{w(\rho,t)}{1 - t}$, the explicit solution is given by
\begin{align*}
    x = \frac{2 (\rho - \frac{1}{2})}{\sqrt{(1 - 2t)^2 + 1}}(1-2t,1)^T ~~\text{and}~~w = (1-t)\rho(\rho-1),
\end{align*}
such that $w$ describes a shrinking parabola and, as in Example 2, $\Gamma(t)$ is a rotating straight line that spans the diameter of $\Omega$. In the left-hand plot of Figure \ref{f:ex2} we display: $\vX^0$ in black, $\vec X^n$, at $t^n = 0.1k$, $k=1,\dots,5$, in blue, and $\partial \Omega$ in red, while in the right-hand plot we display: $W^0$ in black and $W^n$, at $t^n = 0.1k$, $k = 1,\dots,5$ in blue. 
In Table \ref{tab:3} we present the experimental order of convergence for the errors obtained using $\alpha = 0.5$, we do not present the errors for $\alpha=1$ since they are very similar to those obtained using $\alpha=0.5$. For all the four errors we see eocs close to four. 
%where we only consider $\alpha = 0.5$ due to the minute difference between $\alpha = 1$ and $\alpha = 0.5$ in Example 2.  %The computational results are similar to that found in \cite{BDS17,PS15}. 

\begin{remark}
In the three examples presented above if we take $\Dt = C h$  we observe eocs close to two rather than the eocs close to four that we observe above for $\Dt=Ch^2$. Similar convergence behaviour was observed in \cite{BDS17}. 
% as opposed to eocs close to four which we observed for $\Dt = C h^2$. 
\end{remark}

\begin{table}[t]
\caption{Errors and eocs for the parabola defined on the rotating diameter with $\alpha = 0.5$.}%for (\ref{ns})-(\ref{bc2}) .}
\label{tab:3}
    \begin{tabular}{p{0.5cm}p{0.7cm}p{1.4cm}p{0.7cm}p{1.4cm}p{0.7cm}p{1.4cm}p{0.7cm}p{1.4cm}p{0.6cm}}
    \hline\noalign{\smallskip}
    $J$ & $N$ &  $\Er{1}\times e^{-4}$ & $eoc_1$ & $\Er{2}\times e^{-5}$ & $eoc_2$ & $\Er{4} \times e^{-6}$ & $eoc_4$ & $ \Er{5} \times e^{-6} $ & $eoc_5$ \\
    \noalign{\smallskip}\svhline\noalign{\smallskip}
    10 & 50 & 1.205 & - &  3.756 & - & 1.207 & - & 3.073 & - \\
    20 & 200  & 0.07643 & \ter{3.98} & 0.2453 & \ter{3.94} & 0.07829 & \ter{3.95} & 0.2010 & \ter{3.93} \\
    40 & 800 & 0.004795 & \ter{3.99} & 0.01551 & \ter{3.98} & 0.004937 & \ter{3.99} & 0.01271 & \ter{3.98} \\
    80 & 3200 & 0.0003000 & \ter{4.00} & 0.0009721 & \ter{4.00} & 0.0003093 & \ter{4.00} & 0.0007967 & \ter{4.00} \\
    \noalign{\smallskip}\svhline\noalign{\smallskip}
    %160 & 800 & 0.005790 & \ter{2.00} & 0.001208 & \ter{2.00} \\
    %\hline 
    \end{tabular}
\end{table}

\begin{acknowledgement}
JVY gratefully acknowledges the support of the EPSRC grant 1805391. 
VS would like to thank the Isaac Newton Institute for Mathematical Sciences for support and hospitality during the programme {\it Geometry, compatibility and structure preservation in computational differential equations} when work on this paper was undertaken. 
This work was supported by: EPSRC grant number EP/R014604/1.  
\end{acknowledgement}

\end{document}